
\baselineskip=14pt
\parskip=10pt

\magnification=\magstephalf

\def\1{{\overline{1}}}
\def\2{{\overline{2}}}
\parindent=0pt
\overfullrule=0in

\def\frac#1#2{{#1 \over #2}}
\centerline
{\bf 
Approximating Pi (and other constants) by Radicals in the Footsteps 
}
\centerline
{\bf of Rabbi Abraham Ibn Ezra and Guru RSJ Reddy}
\smallskip
\centerline
{\it  Doron ZEILBERGER}
\bigskip
\qquad\qquad\qquad {\it In memory of Gila Aharoni (Dec. 8, 1955-Pi Day, 2025)}
\bigskip
{\bf Abstract:}  Inspired by the polymath Abraham Ibn Ezra (1089-1167)  who approximated $\pi$ by $20/9 \cdot \sqrt{2}$,
and contemporary scholar, RSJ Reddy, who gave $(14-\sqrt{2})/4$, we point out that both of these approximations are optimal in some sense, and we generate many other ones that are even more optimal, albeit not as nice.

\bigskip

{\bf Maple package}: This article is accompanied by a Maple package, {\tt Reddy.txt}, available from

{\tt https://sites.math.rutgers.edu/\~{}zeilberg/tokhniot/Reddy.txt } \quad ,

The front of this article

{\tt https://sites.math.rutgers.edu/\~{}zeilberg/mamarim/mamarimhtml/reddy.html } \quad,

contains numerous output files, some of which will be referred to later.

{\bf Rabbi Abraham Ben Meir Ibn Ezra and Guru RSJ Reddy}

Ibn Ezra was one of the greatest minds of all time, and I already raved about him in [Z]. Since he was an ardent believer in Astrology, he would be considered
a {\it crackpot} by today's mainstream narrow-minded people, but perhaps the jury is still out (in another one thousand years people will realize that
Astrology is valid, when done correctly). At any rate
he is still admired today for his great poetry, philosophy, biblical commentary and  pioneering study of grammar.

In his commentary [I] on Exodus Ch. 3, verse 15, (see also [H] and [Z]),
he gave a construction with ruler and compass  that works out to be
the approximation $\frac{20 \sqrt{2}}{9} $ for $\pi$.

Another great mind, who lives today, and that many people (I am  not one of them) consider a {\it crackpot}, is RSJ Reddy. I am honored to
be in his e-mailing list, and receive, from time to time, yet another approximation to $\pi$ by radicals, and more often, another construction to
an older approximation.
He is what is called, pejoratively, a {\it circle squarer}, but this is wrong! Of course he knows that the
{\it usual} $\pi$ is {\it transcendental}, but he has another, better-looking $\pi$, that he calls the `Cosmic Pi', that equals
$(14-\sqrt{2})/4$. His book [R1], that he very generously allowed me to post on my web-site, contains many such pearls, the most
recent one being [R2].

\vfill\eject

{\bf Approximating Pi the Ibn-Ezra Way}

Ibn-Ezra's approximation has the format
$$
\frac{a}{b} \cdot \sqrt{2} \quad,
$$
where $a$ and $b$ are positive integers.

He gave $\frac{a}{b}=\frac{20}{9}$.

We have

{\bf Prop. 1}: The best approximation of $\pi$ in the form $\frac{a}{b} \sqrt{2}$ with $a$ and $b$ positive integers
and $b<77$ is $\frac{20}{9} \sqrt{2}$.

To see other optimal approximations with the denominator allowed to be larger, see the output file:

{\tt https://sites.math.rutgers.edu/\~{}zeilberg/tokhniot/oReddy4.txt} \quad .

In particular, the convergents of the continued fraction of $\frac{\pi}{\sqrt{2}}$ are among these {\it champions}. Here are the first twenty of them:
$$
2, \frac{9}{4}, \frac{11}{5}, \frac{20}{9}, \frac{311}{140}, \frac{953}{429}, \frac{1264}{569}, \frac{11065}{4981}, \frac{45524}{20493}, \frac{56589}{25474},
\frac{102113}{45967},\frac{158702}{71441}, \frac{260815}{117408},
$$
$$
\frac{1201962}{541073},
\frac{1462777}{658481}, \frac{2664739}{1199554}, \frac{9456994}{4257143}, \frac{12121733}{5456697},\frac{70065659}{31540628},\frac{82187392}{36997325} \quad.
$$

The above output file also has analogous data for approximations of $\pi$ in the form $\frac{a \sqrt{n}}{b}$ for $n \leq 8$, and
readers are welcome to use {\tt Reddy.txt} to generate output for larger $n$ (and other constants).

{\bf Approximating Pi the RSJ Reddy way}

RSJ Reddy's favorite value of $\pi$, that he fondly calls the {\it Cosmic Pi} is $\frac{14-\sqrt{2}}{4}$. It has the format
$$
\frac{a-b\sqrt{2}}{c} \quad,
$$
with $a,b,c$ positive integers. We have

{\bf Prop. 2}: The best approximation of $\pi$ in the form $\frac{a-b\sqrt{2}}{c}$, with $a,b,c$ positive integers, and $b,c\leq 6$ is
Reddy's Cosmic Pi:
$$
\frac{14-\sqrt{2}}{4} \quad.
$$

Procedure {\tt BestRad(C,C1,K)} in our Maple package {\tt Reddy.txt} gives the best Reddy-style approximation of the constant $C$ in terms of
another constant $C1$, in the form $\frac{a-b\,C_1}{c}$ for $1 \leq b,c \leq K$. Here are a few values:

$\bullet$ {\tt BestRad(Pi,sqrt(2),10);} gives $\frac{41}{9}-\sqrt{2}$ with error $0.000250660407332731\dots$ ;

$\bullet$ {\tt BestRad(Pi,sqrt(2),20);} gives $\frac{71}{19}-\frac{8 \sqrt{2}}{19}$ with error $0.00020889037846483\dots$ ;

$\bullet$ {\tt BestRad(Pi,sqrt(2),40);} gives $\frac{52}{11}-\frac{37 \sqrt{2}}{33}$ with error $0.00004668556764564621\dots$ ;

$\bullet$ {\tt BestRad(Pi,sqrt(2),100);} gives $\frac{192}{17}-\frac{98 \sqrt{2}}{17}$ with error $ 2.484464588\dots \times 10^{-7}\dots$ .

{\bf Using the PSLQ algorithm to generate Reddy-style approximations to Pi (or any other constant)}

Procedure {\tt BestRad} is {\it brute-force}, trying out all the possible values of $1 \leq b,c \leq K$. Of course it can be optimized,
but a much more efficient way to generate, very fast, Reddy-style approximations is the amazing {\tt PSLQ} algorithm [FBA] invented
by experimental mathematician David H. Bailey and  mathematician-sculptor Helaman Ferguson (see also  [BB]). PSLQ is implemented
very efficiently in Maple, and the fast version of {\tt BestRad}, called {\tt BestRadF(C,C1,d)}, asks PSLQ to find an integer relation
between (floating-point approximations of) $1,C,C1$. Alas, we can not target the best such approximation with $b,c \leq K$, like we did before, but have to be happy
with what we get. In order to get nice-looking approximations we get find the relation between $1,C,C1$ with the {\tt Digits} parameter
set smaller than the default $10$ (or larger). This is the parameter {\tt d} in {\tt BestRadF}.

Typing

{\tt seq(BestRadF(Pi,sqrt(2),i),i=2..10);}

yields
$$
1+\frac{3 \sqrt{2}}{2}, -\frac{5}{2}+4 \sqrt{2}, -\frac{74}{3}+\frac{59 \sqrt{2}}{3}, \frac{49}{24}+\frac{7 \sqrt{2}}{9}, \frac{123}{176}+\frac{19 \sqrt{2}}{11}, 
-\frac{227}{88}+\frac{89 \sqrt{2}}{22},
$$
$$
\frac{9509}{3856}+\frac{921 \sqrt{2}}{1928}, \frac{8649}{2576}-\frac{295 \sqrt{2}}{1932}, 
\frac{95154}{35671}+\frac{11957 \sqrt{2}}{35671} \quad .
$$

Readers are welcome to find many other Reddy-style approximations not just of $\pi$ and not just in terms of $\sqrt{2}$ or $\sqrt{n}$.

For example, typing

{\tt BestRadF(Zeta(3),gamma,10);}

will tell you that $\zeta(3)$ is approximately
$$
\frac{84856}{13637}-\frac{118610 \gamma}{13637} \quad,
$$
where $\gamma$ is Euler's  constant. This agrees up to the tenth decimal place.

If you want agreement up to twenty places, change the $10$ to $20$ and get:
$$
\frac{5503504394}{8892617755}+\frac{1796877140 \gamma}{1778523551} \quad .
$$

{\bf Conclusion}: Inspired by two of my heroes, I wrote a Maple package {\tt Reddy.txt}, that can reproduce their nice approximations to $\pi$,
prove optimality, in some sense, and generate many other approximations, much more optimal, in other senses.


{\bf References}

[BB] David H. Bailey and Jonathan M. Borwein, {\it PSLQ: An algorithm to discover integer relations} \hfill\break
{\tt https://www.davidhbailey.com/dhbpapers/pslq-comp-alg.pdf} \quad .

[FBA] Helaman R. P. Ferguson, David H. Bailey and Stephen Arno, {\it Analysis of PSLQ, An Integer Relation Finding Algorithm},
Math. of Computation {\bf 68} (1999), 351-369.

[H] Haim Hanani, {\it Letter to the Editor}, Gilyonot LeMatematika [edited by Joseph Gillis], Vol. {\bf 4}, No. {\bf 1} (Nov. 1969). [in Hebrew] \hfill\break
[Here is a scan of the letter:{\tt  http://sites.math.rutgers.edu/\~{}zeilberg/gilyonot/hanani.pdf}, and here is a scan of the full issue: \hfill\break
{\tt https://sites.math.rutgers.edu/\~{}zeilberg/gilyonot/gilyonotNov1969.pdf}.

[I] Abraham Ibn Ezra, {\it Commentary on Exodus Ch. 3, verse 15} [in Hebrew]: (given, for example, in ``Mikraot Gdolot'', sefer Shmot, rav Pninim,
Levin-Epstein, Jerusalem.

[R1] R. Dhanalakshmi Sarva Jagannadha Reddy, {\it ``PI OF THE CIRCLE''}, 2014 \hfill\break
available here: {\tt https://sites.math.rutgers.edu/\~{}zeilberg/akherim/Reddy2014.pdf} \quad.

[R2] R. Dhanalakshmi Sarva Jagannadha Reddy, {\it ``ARE CIRCUMFERENCE AND ITS Pi NUMBER LOCKED IN
RADIUS JUST LIKE ENERGY LOCKED INSIDE MASS?''}, 28 Aug. 2026 \hfill\break
available here: {\tt https://sites.math.rutgers.edu/{}\~{}zeilberg/akherim/reddy26.pdf} \quad.

[Z] Doron Zeilberger, {\it King Solomon and Rabbi Ben Ezra's Evaluations of Pi and Patriarch Abraham's Analysis of an Algorithm} \hfill\break
{\tt https://sites.math.rutgers.edu/\~{}zeilberg/mamarim/mamarimhtml/king.html} \quad .

\bigskip
\hrule
\bigskip
 Doron Zeilberger, Department of Mathematics, Rutgers University (New Brunswick), Hill Center-Busch Campus, 110 Frelinghuysen
Rd., Piscataway, NJ 08854-8019, USA. {\tt DoronZeil@gmail.com} \hfill\break

{\bf Sept. 9, 2026} 

\end